\documentclass[conference]{IEEEtran}
\IEEEoverridecommandlockouts
\usepackage{cite}
\usepackage{amsmath,amssymb,amsfonts}
\usepackage{siunitx}
\usepackage{algorithmic}
\usepackage{graphicx}
\usepackage{textcomp}
\usepackage{xcolor}
\usepackage{booktabs}
\usepackage{bm}
\def\BibTeX{{\rm B\kern-.05em{\sc i\kern-.025em b}\kern-.08em
    T\kern-.1667em\lower.7ex\hbox{E}\kern-.125emX}}
\usepackage{xurl}
\usepackage[ruled]{algorithm2e}
\SetAlCapNameFnt{\footnotesize}
\SetAlCapFnt{\footnotesize}
\SetAlFnt{\small}

\usepackage[acronym]{glossaries}
\newacronym{scopf}{SCOPF}{security constrained optimal power flow}
\newacronym{rmac}{RMAC}{reliability management approach and criterion}
\newacronym{smw}{SMW}{Sherman-Morrison-Woodbury}
\newacronym{drc}{DRC}{deterministic reliability criterion}
\newacronym{minlp}{MINLP}{mixed-integer nonlinear program}
\newacronym{tso}{TSO}{transmission system operator}
\newacronym{imml}{IMML}{inverse matrix modification lemma}
\newacronym{ptdf}{PTDF}{power transfer distribution factors}
\newacronym{lodf}{LODF}{line outage distribution factors}

\begin{document}

\title{
A Novel Approach for Solving Security Constrained Optimal Power Flow Using the Inverse Matrix Modification Lemma and Benders Decomposition
\thanks{The research leading to results in this paper has received funding from the Research
Council of Norway through the project “Resilient and Probabilistic reliability
management of the transmission grid” (RaPid) (Grant No. 294754), The
Norwegian Water Resources and Energy Directorate, and Statnett.}
}

\author{\IEEEauthorblockN{Matias Vistnes and Vijay Venu Vadlamudi}
\IEEEauthorblockA{ \textit{Department of Electric Energy} \\
\textit{Norwegian University of Science and Technology}, \\
Trondheim, Norway \\
\{matias.vistnes, vijay.vadlamudi\} @ntnu.no}
\and
\IEEEauthorblockN{Sigurd Hofsmo Jakobsen and Oddbjørn Gjerde}
\IEEEauthorblockA{
\textit{SINTEF Energy Research}\\
Trondheim, Norway \\
\{sigurd.jakobsen, oddbjorn.gjerde\} @sintef.no}
}

\maketitle

\begin{abstract}
With the increasing complexity of power systems, faster methods for power system reliability analysis are needed.
We propose a novel methodology to solve the security constrained optimal power flow (SCOPF) problem that reduces the computational time by using the Sherman-Morrison-Woodbury identity and Benders decomposition. The case study suggests that in a 500 node system, the run time is reduced by 83.5\% while ensuring a reliable operation of the system considering short- and long-term post-contingency limits and reducing the operational costs, compared to a preventive `N-1' strategy.
\end{abstract}


\section{Introduction}
Even with increased processing capability of computers, all the desired modeling details in power system analysis cannot be realized without computational bottlenecks in the ever-growing power systems.
More concerning for system operators and planners is that the system's state is changing faster as a result of higher penetration of renewable energy and more frequent extreme weather.
Also, with the energy transition and electrification, power system investment is not keeping up with demand growth, leading to operation closer to the equipment limits.
Consequently, it is vital to have accurate solutions that can ensure the reliability of electric supply for customers; society is more dependent on it than ever before.

Using the traditional approach of reliability management, reliability could be accomplished through grid development solutions that meet a \gls{drc} at all times (typically \mbox{N-1}\footnote{``The N-1 criterion is a principle according to which the system should be able to withstand at all times a credible contingency - i.e., unexpected failure or outage of a system component (such as a line, transformer, or generator) - in such a way that the system is capable of accommodating the new operational situation without violating operational security limits'' \cite{garpur_generally_2013}.}). Preventive rescheduling of the economic dispatch is often used by operators to meet a \gls{drc}. However, such solutions may be prohibitively expensive and not necessarily socio-economic efficient. Recommendations toward more probabilistic and risk-based methods to reliability management, balancing reliability and costs, were made in the EU-funded project GARPUR \cite{garpur_generally_2013}, where the proposed \gls{rmac} contains a socio-economic objective, a reliability target, a discarding principle, and a relaxation principle. 
Recent research shows the need for better methods for \gls{rmac}, especially the need for methods that scale better with the system size \cite{garpur_generally_2013, capitanescu_critical_2016, kardos2020schur}. 

\Gls{scopf} is an extensively used tool in \gls{rmac} studies. It can be used to optimize the generation schedule and power system operation. A comprehensive \gls{scopf} formulation should include preventive, corrective, and restorative actions to meet system limits both in the current state and in the potential contingency states using an AC power system model. It is a large-scale non-convex \gls{minlp} \cite{capitanescu_critical_2016}, and no all-inclusive solution methods exist to analyze large realistic systems \cite{Wang2016, karangelos_iterative_2019}. There are many sources of computational burden which make the problem intractable; especially system size, number of contingencies considered, and the optimization of actions---preventive, corrective, and restorative. Risk and chance constraints increase the complexity of the problem \cite{Wang2016}. Capitanescu \cite{capitanescu_critical_2016} noted that several sets of post-contingency limits (on different timescales) should be included to better model the actions taken by operators and the different timescales of equipment operating limits. 
Many \gls{scopf} formulations are proposed in the literature. Wang \textit{et al.} \cite{Wang2016} solve a risk based DC \gls{scopf} for a large power system using Lagrangian relaxation and Benders decomposition. Kardos \textit{et al.} \cite{kardos2020schur} use parallel, distributed algorithms to solve a non-probabilistic AC \gls{scopf} with preventive actions using Schur-Complements on a large scale system. Karangelos and Wehenkel \cite{karangelos_iterative_2019} solve a chance-contrained AC \gls{scopf} including preventive and corrective actions on a smaller power system using an iterative scheme.

Solving an \gls{scopf} for large power systems involves a large number of nodes, branches, and contingencies considered. Two approaches are focused here: the \gls{imml} and Benders decomposition. The \gls{imml} shows how to efficiently compute matrix small modifications in large matrices, first known as the Sherman-Morrison-Woodbury identity \cite{Woodbury1950}, and later used in the power system context \cite{alsac1983sparsity}. 
While underused in literature, the \gls{imml} is an application of Schur-Complements that could have great potential for use in \gls{scopf} formulations through its reduced calculation burden of contingency cases. The \Gls{imml} is combined with Benders decomposition, a widely used, highly efficient method to decompose large optimization problems \cite{Benders1962, geoffrion1972, Wang2016}. 

This paper proposes a new, efficient and scalable methodology to find the optimal socio-economic operation of a power system using both preventive and corrective actions. Restorative actions are planned to be incorporated in future work. The proposal integrates the \gls{imml} and Benders decomposition into a probabilistic DC \gls{scopf} framework; the developed method is modular and can be combined with other \gls{scopf} formulations and solution techniques. A DC \gls{scopf} is used to enable analysis of large power systems. The central contribution lies in the combination of efficient solving of the power flow and the unique formation of Benders' cuts including both preventive and corrective actions.

The rest of the paper is organized as follows: Section \ref{sec:method} presents the \gls{scopf} formulation, the \gls{imml}, and Benders decomposition, resulting in the proposed methodology. Section \ref{sec:case_study} presents the details of the case study. Section \ref{sec:conclusion} presents a concluding summary of the work.

\section{Methodology}\label{sec:method}
The power system consists of sets of nodes $\mathcal{N}$, branches $\mathcal{B}$, generators $\mathcal{G}$, and demands $\mathcal{D}$, where $N = |\mathcal{N}|$ and $B = |\mathcal{B}|$. Contingencies are denoted $c$ and taken from the contingency set $\mathcal{C}$. 
To capture short- and long-term operational limits on transmission lines and the effect of generator ramping, we consider a power system with one pre-contingency state and two post-contingency states. The first post-contingency state is set after the contingency when circuit breakers have tripped. The second state is set minutes later when all frequency reserves are activated. We assume that the system has a stable trajectory between all states; although this assumption is common, it is not always valid \cite{BODAL2022Probabilistic}. 

\subsection{Security Constrained Optimal Power Flow (SCOPF)}\label{ssec:scopf}
A preventive-corrective \gls{scopf} from Capitanescu and Wehenkel
\cite{Capitanescu2007improving} is extended, and presented below (matrix and vector variables are in bold):
\begin{subequations}\label{eq:scopf}
\begin{align}
    \min_{\bm{x}_0,\bm{u}_0,\bm{x}_c,\bm{u}_c} \quad & f = f(\bm{x}_0,\bm{u}_0,\bm{x}_c,\bm{u}_c)   \label{eq:scopf_f}\\
    s.t. \quad & \bm{g_0}(\bm{x}_0,\bm{u}_0) = 0  \label{eq:scopf_g0}\\
    & \bm{h_0}(\bm{x}_0,\bm{u}_0) \leq \bm{\overline h}_0^{LT}  \label{eq:scopf_h0}\\
    & \bm{g_c}(\bm{x}_c,\bm{u}_c) = 0 &\quad c \in \mathcal{C}  \label{eq:scopf_gk}\\
    & \bm{h_c}(\bm{x}_c,\bm{u}_0) \leq \bm{\overline h}_c^{ST} &\quad c \in \mathcal{C}   \label{eq:scopf_hk}\\
    & \bm{h_c}(\bm{x}_c,\bm{u}_c) \leq \bm{\overline h}_c^{LT} &\quad c \in \mathcal{C} \\
    & |\bm{u}_c - \bm{u}_0| \leq \Delta \bm{u}_c &\quad c \in \mathcal{C} \label{eq:scopf_delta}
\end{align}
\end{subequations}
where $\bm{x}$ is the state variables, $\bm{u}$ is the control variables, $f$ is the function to be minimized with respect to $\bm{x}_0$, $\bm{u}_0$, $\bm{x}_c$, and $\bm{u}_c$, $\bm{g}$ is the power flow equations, $\bm{h}$ is the operating limits, and $\Delta \bm{u}_c$ is the maximum allowed change of control variables.
Subscripts $\cdot_0$ and $\cdot_c$, respectively, represent base case and post-contingency case. $\bm{\overline h}_0^{LT}$, $\bm{\overline h}_c^{ST}$, and $\bm{\overline h}_c^{LT}$ are operating limits for normal operation, short-term limits after a contingency, and long-term limits after a contingency, respectively.

In this paper branch contingencies are considered, which include transmission lines and transformers. The operational limits considered are maximum current in branches (short- and long-term), generator maximum active power output, generator ramping limit, and a maximum of 10\% load shedding on each node after a contingency. 

\subsection{Inverse matrix modification lemma}\label{ssec:imml}
When performing a contingency analysis we solve the DC power flow equation: 
\begin{equation}\label{eq:base}
    (\bm{H} + \bm{\Delta H}) \cdot \bm{\theta} = \bm{P}
\end{equation}
where $\bm{H}$ is a sparse ($N\times N$) susceptance matrix, $\bm{\Delta H}$ is a system modification due to a contingency, $\bm{\theta}$ is the node voltages angles, and $\bm{P}$ is the node injected active power. 
If $\bm{\Delta H}$ is symmetric\footnote{For an asymmetric $\bm{\Delta H}$ (e.g., when using phase shifting transformers), the formula is still valid but the equations differ slightly \cite{alsac1983sparsity}.}, and can be written as:
\begin{equation}
    \bm{\Delta H} = \bm{\Phi}_m \cdot \bm{\delta h} \cdot \bm{\Phi}_m^\mathsf{T}
\end{equation}
where $\bm{\delta h}$ is an ($M \times M$) matrix of the modifications in $\bm{H}$, $\bm{\Phi}_m$ is an ($N \times M$) connectivity matrix, and $\bm{\Phi}^\mathsf{T}$ is the transposed matrix of $\bm{\Phi}$. $M$ is the number of modified components. 
This is a rank-one updated problem that allows us to use the \gls{imml}. The \gls{imml} facilitates a cheaper computation than solving the full equation system. Further, solving for $\bm{\theta}$ and using the \gls{imml} yields \cite{alsac1983sparsity}:
\begin{gather}\label{eq:imml}
    \bm{\theta} = \left(\bm{I} - \bm{H}^{-1} \bm{\Phi}_m \bm{c} \bm{\Phi}_m^\mathsf{T}\right) \bm{H}^{-1} \bm{P} \\
    \bm{c} = \left(\bm{\delta h}^{-1} + \bm{\Phi}_m^\mathsf{T} \bm{H}^{-1} \bm{\Phi}_m\right)^{-1}
\end{gather}
where $\bm{I}$ is the identity matrix. \eqref{eq:imml} is valid for any number of system element modifications. In the case of a single component modification, $\bm{c}$ becomes a scalar, $c$. In general, if $\bm{c}^{-1}$ is singular the system modification separates the system into islands \cite{alsac1983sparsity}. Separated systems cannot be analyzed using the \gls{imml}.
In the Big O notation, the scheme is O($NM+M^2$), greatly reducing the number of operations when $N \gg M$, compared to solving the linear problem, O($N^3$).

The procedure using the \gls{imml} to calculate the branch power flow after contingency $\bm{F}_c$ of branch $l$ from node $i$ to node $j$ contains the following equations:
\begin{gather}
    \bm{\theta}_{c} = \bm{\theta} - 
    \frac{\bm{\delta} \cdot \left(\theta[i] - \theta[j]\right)}{1/H[i,j] + \delta[i] - \delta[j]} \\
    \bm{\delta} = \frac{x_l}{H[i,j]} \cdot \left(\bm{X}[:,i] - \bm{X}[:,j]\right) \\
    \bm{F}_{c} = \bm{\Psi} \cdot \bm{\Phi} \cdot \bm{\theta}_{c} \label{eq:vm_new}
\end{gather}
where $x_l$ is the reactance of branch $l$, $\bm{X}$ is the inverse susceptance matrix, $\bm{\Psi}$ is the diagonal branch susceptance matrix, and \mbox{$\bm{X}[:,i]$} ($\bm{X}[i,:]$) is notation for the $i$th column (row) of matrix $\bm{X}$. Details on the matrices are given in the appendix, Section~\ref{sec:Appendix}. 

In the proposed methodology, the inverse susceptance matrix after a contingency, $\bm{X}_{c}$, is used to calculate contingency \gls{ptdf} matrix (effectively the \gls{lodf}, $\bm{\varphi}$) and the power flow after the outage.
\begin{gather}
    \bm{\varphi}_{c} = \bm{\Psi} \cdot \bm{\Phi} \cdot \bm{X}_c \\
    \bm{F} = \bm{\varphi}_{c} \cdot \bm{P} \label{eq:phi_new}
\end{gather}
Direct calculation of $\bm{X}_{c}$ is done by inverting the contingency susceptance matrix, $\bm{X}_{c} = (\bm{H} + \bm{\Delta H})^{-1}$. However, using the \gls{imml} to find $\bm{X}_{c}$ is more efficient.
\begin{gather}
    \bm{X}_{c} = \bm{X} - \frac{\bm{X}\cdot \bm{\Phi}[i,:]\cdot H[i,j]\cdot \bm{\delta}}{1+H[i,j]\cdot \bm{\delta}\cdot \bm{\Phi}[i,:]} \\
    \bm{\delta} = \frac{x_b}{H[i,j]} \cdot \left(\bm{X}[:,i] - \bm{X}[:,j]\right)
\end{gather}

\subsection{Benders decomposition}\label{ssec:benders}
Benders decomposition is a method to divide an optimization problem into a main problem and several sub-problems \cite{Benders1962}. From the sub-problems a Benders' cut is added to the main problem formulation in such a way that an optimal solution of the main problem is a feasible solution of \eqref{eq:scopf}.
In our case the main problem is \eqref{eq:scopf_f} -- \eqref{eq:scopf_h0} and the constraints in the sub-problem are \eqref{eq:scopf_gk} -- \eqref{eq:scopf_delta}. Often the contingencies with active constraints in the optimal solution are a small sub-set of $\mathcal{C}$. Thus not all sub-problems need to be added to the main problem to make a feasible solution.
The branch power flow constraints in the sub-problem can be reformulated to reduce the complexity through a Benders' cut.

Starting from the equation for branch power flow \eqref{eq:phi_new}, the Benders' cut can be deduced.
Using \eqref{eq:phi_new} a constraint can be set up to limit the branch power flow to the rated value $\bm{\overline h_l}$: 
\begin{equation}
    -\bm{\overline h_l} \leq \bm{\varphi} \cdot \bm{P} \leq \bm{\overline h_l}.
\end{equation}
Only the lower or the upper bound can be active in the optimal solution. \eqref{eq:phi_new} is also valid for a change in injected power $\Delta P$. Thus, a constraint which requires a power flow change $\Delta F$ can be expressed.
\begin{gather}
    \bm{\varphi} \cdot \bm{\Delta P} \ \left\{ 
    \begin{array}{@{}r@{\thinspace}l}
        \geq \Delta F_l ,&\quad\text{if } F_l \geq 0 \\
        \leq \Delta F_l ,&\quad\text{if } F_l < 0
    \end{array}
    \right. \quad \forall l \in \mathcal{B}
\end{gather}
where $F_l$ is the current flow on branch $l$.
The Benders' cut added to the main problem formulation in \eqref{eq:scopf_f} -- \eqref{eq:scopf_h0} to mitigate overload on branch $l$ after a contingency on branch $i$ is given as follows:
\begin{gather} \label{eq:benders_cut}
    \sum^{\mathcal{N}}_{n} \varphi_{i,n} \Delta P_{n} \ \left\{
    \begin{array}{@{}r@{\thinspace}l}
        \geq F_l - \overline h_l ,&\quad\text{if } F_l > \overline h_l \\
        \leq F_l + \overline h_l ,&\quad\text{if } F_l < - \overline h_l
    \end{array}
    \right. \\ \label{eq:DeltaF}
    \Delta P_{n} = P_{c,n} - P_{0,n}^g - P_{c,n}^{g+} +P_{c,n}^{g-} + P_{c,n}^d
\end{gather}
where $\varphi_{i,n}$ is the value of the \gls{lodf}
at the contingency branch $i$ and node $n$, $P_{c,n}$ is the injected active power at node $n$ using the current generation schedule solution for contingency $c$, $P_{0,n}^g$ and $P_{c,n}^{g}$ are the active power generation change at node $n$ as a preventive action and as a corrective action, respectively ($P^{g+}$ is an increase and $P^{g-}$ is a decrease), and $F_l$ is the flow on branch $l$ with overload. \eqref{eq:benders_cut} is not defined for branches within their limit, $- \overline h_l \leq F_l \leq \overline h_l$, as there is no overload to mitigate.

If only considering preventive actions or only short-term operating limits, \eqref{eq:DeltaF} is reduced to, respectively, as follows:
\begin{gather}
    \Delta P_{n} = P_{0,n} - P_{0,n}^g \\
    \Delta P_{n} = P_{c,n} - P_{0,n}^g + P_{c,n}^{g-} + P_{c,n}^d
\end{gather}
where $P_{0,n}$ is the injected active power at node $n$ using the current generation schedule solution in the base case. The corresponding $\overline h_l$ for the time-frame is used.

In addition to the Benders' cut modelling the flow of power, the objective function is extended using \eqref{eq:obj_corrective} to account for the added costs of mitigating the branch overload when using corrective actions, or \eqref{eq:obj_short-term} for only short-term operating limits.
\begin{gather}\label{eq:obj_corrective}
    \min k_{\mathcal{G}}^\mathsf{T} \left(P^{\mathcal{G}}_0 +  \pi_{c} P^{\mathcal{G}+}_{c}\right) + \text{voll}^\mathsf{T} \left(P^\mathcal{D}_0 + \pi_{c} P^\mathcal{D}_c\right)\\ \label{eq:obj_short-term}
    \min \text{voll}^\mathsf{T} \left(P^\mathcal{D}_0 + \pi_{c} P^\mathcal{D}_c\right)
\end{gather}
where $k^{\mathcal{G}}$ is the generation cost, $\pi_{c}$ is the contingency probability, voll is the value-of-lost-load, and $P^{\mathcal{D}}$ is the active power load shedding.

\subsection{Algorithm}
\begin{figure}[tb]
    \centering
    \includegraphics[width=\linewidth]{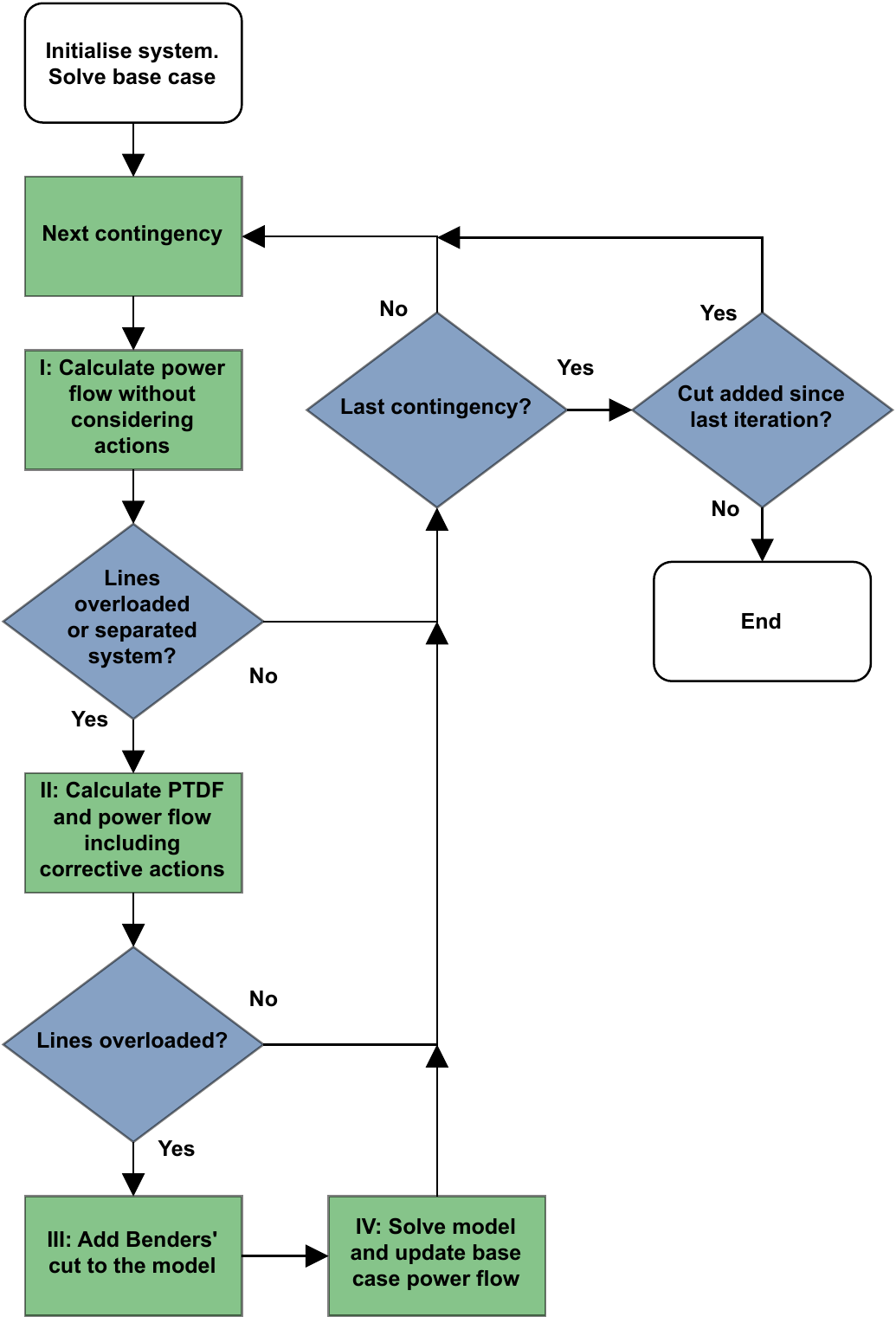}
    \caption{The proposed algorithm for solving \gls{scopf}.}
    \label{fig:SCOPF}
\end{figure}

The proposed methodology, illustrated in Fig.~\ref{fig:SCOPF}, employs the \gls{imml} to identify contingencies where branch overloads may occur. Benders' cuts are added to the model to represent potential preventive or corrective actions that can be taken to mitigate overloaded branches.

Initially, an \gls{scopf} is solved without considering any possible contingencies, using any algorithm that optimizes linear programs. This generates an `N-0' safe generation schedule for the system at its current state, and gives the lowest system cost. As cuts are later added as constraints, the system cost increases.

For a given branch contingency, the power flow is calculated with the base case generation schedule (with possible preventive actions) using \eqref{eq:vm_new}. If the calculation finds overloaded branches from the contingency, \eqref{eq:phi_new} is calculated. \eqref{eq:phi_new} is O($n^2m+nm^2$) while \eqref{eq:vm_new} is considerably faster at O($nm+m^2$).

If the system is separated into islands after a branch contingency, Goderya's algorithm \cite{goderya1980fast} is used to detect them. Then the reference island (the connected nodes which includes the reference node) is found. Islands without the reference node are assumed to black out. The reference island nodes of $\bm{H}$ are used in Algorithm~\ref{alg:X} to calculate $\bm{X}_c$; then \eqref{eq:phi_norm} and \eqref{eq:flow_norm} are used to find the new \gls{ptdf} matrix and power flow, respectively, followed by the Benders' cut if there are any overloaded branches. If not, short-term variables are added to the main problem, though no cut is added to the formulation.

When a contingency results in the overload of one or more branches, Benders' cuts for preventive and corrective actions are added for each overload, as outlined in section \ref{ssec:benders}. The model is then resolved to find the new optimal point and update the base case power flow.

The algorithm ends once all contingencies have been checked, and no overloaded branches are found. The output is an optimal generation schedule for a base case together with optimal corrective actions after a contingency for the optimal system cost.

\section{Case study}\label{sec:case_study}
A case study is run to compare a standard \gls{scopf} (direct optimization of \eqref{eq:scopf}) and the proposed methodology. The methods are implemented in Julia language \cite{Bezanson2017Julia} with the packages JuMP \cite{Lubin2023JuMP} and PowerSystems \cite{Lara2021PowerSystems} and problem optimization using the Gurobi solver \cite{gurobi}.
Results and run time are compared for the IEEE RTS-79 system \cite{Subcommittee1979rts} and the SouthCarolina500 (ACTIVSg500) test system \cite{Birchfield2017g2000} in Table~\ref{tab:runtime}. Where generator ramping limits were not provided, limits of 1\% of rated power per minute were set. The results show that using the proposed methodology significantly reduces the solution time while yielding highly similar optimal cost values. 

\begin{table}[tb]
    \centering
    \caption{Optimal cost value and run time using the BenchmarkTools package in Julia. Run on \textregistered Intel \texttrademark Core i5-1145G7.}
    \begin{tabular}{@{}lrrrr@{}} \toprule
         & \multicolumn{2}{c}{Standard Methodology} & \multicolumn{2}{c}{Proposed Methodology} \\ 
         & Optimal Cost & Time [s] & Optimal Cost & Time [s] \\ \cmidrule(l){2-5} 
        IEEE RTS-79 & 483.00 & 0.106 & 483.00 & 0.022  \\
        ACTIVSg500 & 3243 & 298.091 & 998.64 & 49.154 \\
        \bottomrule
    \end{tabular}
    \label{tab:runtime}
\end{table}

One hypothesis in the proposed methodology is that it is faster to calculate  only the \gls{ptdf} matrix when needed, instead of calculating it for every contingency. For the \mbox{ACTIVSg500} system, calculating contingency branch power flow using $\bm{\theta}_c$ in \eqref{eq:vm_new} takes \SI{2.6}{\micro\second}, while it is much slower using $\bm{\varphi}_c$ in \eqref{eq:phi_new} (or \eqref{eq:phi_norm}), taking \SI{1.8}{\milli\second} (\SI{7.0}{\milli\second}). When applying the proposed methodology on the ACTIVSg500 system, approximately 1/3rd of the run time is used by the solver to optimize the model and approximately 1/5th of the run time is used to compute the branch power flow and \gls{ptdf}. If skipping \eqref{eq:vm_new}, and always calculating the branch power flow using $\bm{\varphi}_c$, the branch power flow run time time increases. When solving the RTS-79 system using the proposed methodology, 1/4th of the run time is the solver run time and 1/10th of the run time is the branch power flow run time; branch power flow run time increases when skipping \eqref{eq:vm_new}. The RTS-79 system, which has 38 branches, has only one single branch contingency case that results in the system separation. The ACTIVSg500 system, which has 597 branches, has 268 single branch contingency cases resulting in the system separation. Thus, even on systems with a high ratio of branch contingencies that result in separation of the system, first calculating $\bm{\theta}_c$ and only calculating $\bm{\varphi}_c$ when needed is shown to decrease the run time.

The IEEE RTS-79 system is used to compare the proposed and the standard methodology. For the proposed methodology, the Benders' cuts of branch contingencies are shown in Table~\ref{tab:cuts}. Branches are denoted as follows: branch from node $i$ to node $j$ is i-j.

\begin{table}[tb]
\centering
\caption{Benders' cuts added while solving the IEEE RTS-79 system.}
\label{tab:cuts}
\begin{tabular}{@{}llrl@{}}
\multicolumn{4}{c}{\textbf{Iteration 1}} \\ \toprule
Type & Contingency & \multicolumn{2}{l}{Overload @ branch [pu]} \\ \midrule
Corrective & 3-24 & 0.7391 &@ 14-16 \\
Corrective & 12-23 & 0.2420 &@ 14-16 \\
Corrective & 13-23 & 0.3652 &@ 14-16 \\
Corrective & 14-16 & 0.0491 &@ 3-24 \\
Corrective & 15-21 & 0.1415 &@ 16-17 \\
Corrective & 15-21 & 0.1415 &@ 16-17 \\
Corrective & 15-24 & 0.7391 &@ 14-16 \\
Corrective & 16-19 & 0.2183 &@ 14-16 \\ \midrule
\\
\multicolumn{4}{c}{\textbf{Iteration 2}} \\ \toprule
Type & Contingency & \multicolumn{2}{l}{Overload @ branch [pu]} \\ \midrule
Preventive & 3-24 & 0.0700 &@ 7-8 \\
Corrective & 3-24 & 0.3500 &@ 7-8 \\
Preventive & 12-23 & 0.0700 &@ 7-8 \\
Corrective & 12-23 & 0.3500 &@ 7-8 \\
Preventive & 13-23 & 0.0700 &@ 7-8 \\
Corrective & 13-23 & 0.3500 &@ 7-8 \\
Corrective & 14-16 & 0.1261 &@  7-8 \\
Corrective & 15-21 & 0.2234 &@ 7-8 \\
Corrective & 15-21 & 0.2234 &@ 7-8 \\
Preventive & 15-24 & 0.0700 &@ 7-8 \\
Corrective & 16-19 & 0.3217 &@ 7-8 \\ \bottomrule
\end{tabular}
\end{table}
From Table~\ref{tab:cuts}, it is evident from the first iteration of the algorithm that all overloads on the initial base case generation schedule are mitigated, for all contingencies. The second iteration shows overload only on branch \mbox{7-8}. These overloads are generated by corrective actions to mitigate the overloads found in iteration one. By analyzing the Benders' cut \eqref{eq:benders_cut} it can be seen that a power injection change at any node can be used to fulfill the constraint. However, only the rating of the overloaded branch is used to make the cut, and thus other branches can be overloaded by the corrective actions found after adding the Benders' cut. In theory, the algorithm would need $B^2$ iterations to ensure that all branch ratings are within limits. In the study of both the RTS-79 system and the ACTIVSg500 system, this resulted in one additional iteration, two iterations in total for each. 

A comparison of the generation schedule produced by the proposed and the standard methodology is shown in Table~\ref{tab:P_contingency}. The schedule is optimized on a per generator and per load basis, but aggregated to per node basis in the tables for the sake of compactness. The base case shows equal injected power between the methods for all nodes (not shown). Both methods find optimal short-term corrective actions for the contingency of branch 7-8. This contingency will isolate node 7 from the rest of the system, and thus actions must be taken to account for the lost generation on node 7. While the generation schedule is different, the optimal cost is equal down to the fifth significant digit.

\begin{table}[tb]
\centering
\caption{Injected power changes per node after contingencies.}
\label{tab:P_contingency}
\begin{tabular}{@{}llr@{}l@{}llr@{}}
\multicolumn{7}{c}{\textbf{Short-term Corrective Actions}}    \\ \toprule
\multicolumn{3}{c}{Standard Methodology}  & \phantom{$\qquad$}  & \multicolumn{3}{c}{Proposed Methodology} \\
Contingency   & Node  & Change &  & Contingency   & Node  & Change \\ \cmidrule{1-3} \cmidrule{5-7}
7-8   & 3    & 0.18   &  & 3-24   & 3    & 0.07   \\
7-8   & 4    & 0.07   &  & 3-24   & 7    & -0.07  \\
7-8   & 7    & -1.40   &  & 7-8   & 3    & 0.18   \\
7-8   & 8    & 0.17   &  & 7-8   & 4    & 0.07   \\
7-8   & 10   & 0.20    &  & 7-8   & 7    & -1.40   \\
7-8   & 13   & 0.27   &  & 7-8   & 8    & 0.17   \\
7-8   & 15   & 0.32   &  & 7-8   & 10   & 0.20    \\
7-8   & 16   & 0.10    &  & 7-8   & 13   & 0.27   \\
7-8   & 20   & 0.10    &  & 7-8   & 15   & 0.32   \\ \cline{1-3}
            &      &        &  & 7-8   & 16   & 0.10    \\
            &      &        &  & 7-8   & 20   & 0.10    \\
            &      &        &  & 12-23 & 3    & 0.07   \\
            &      &        &  & 12-23 & 7    & -0.07  \\
            &      &        &  & 13-23 & 3    & 0.07   \\
            &      &        &  & 13-23 & 7    & -0.07  \\
            &      &        &  & 15-24 & 3    & 0.07   \\
            &      &        &  & 15-24 & 7    & -0.07  \\
            &      &        &  & 16-19 & 3    & 0.04   \\
            &      &        &  & 16-19 & 7    & -0.04  \\ \cline{5-7}
\\
\multicolumn{7}{c}{\textbf{Long-term Corrective Actions}}     \\ \toprule
\multicolumn{3}{c}{Standard Methodology}  &  & \multicolumn{3}{c}{Proposed Methodology} \\  
Contingency   & Node  & Change &  & Contingency   & Node  & Change \\ \cmidrule{1-3} \cmidrule{5-7}
3-24   & 2    & 0.39   &  & 3-24   & 2    & 0.39   \\
3-24   & 3    & 0.18   &  & 3-24   & 3    & 0.18   \\
3-24   & 13   & 0.90    &  & 3-24   & 13   & 0.90    \\
3-24   & 18   & -1.47  &  & 3-24   & 18   & -1.47  \\
7-8   & 1    & 0.20    &  & 7-8   & 7    & -1.40   \\
7-8   & 2    & 0.20    &  & 7-8   & 13   & 0.90    \\
7-8   & 3    & 0.18   &  & 7-8   & 15   & 0.50    \\
7-8   & 7    & -1.40   &  & 12-23 & 13   & 0.55   \\
7-8   & 8    & 0.02   &  & 12-23 & 16   & -0.30   \\
7-8   & 13   & 0.30    &  & 12-23 & 18   & -0.25  \\
7-8   & 15   & 0.50    &  & 13-23 & 13   & 0.71   \\
12-23 & 13   & 0.55   &  & 13-23 & 16   & -0.30   \\
12-23 & 16   & -0.30   &  & 13-23 & 18   & -0.41  \\
12-23 & 18   & -0.25  &  & 14-16 & 13   & 0.16   \\
13-23 & 13   & 0.71   &  & 14-16 & 15   & -0.16  \\
13-23 & 16   & -0.30   &  & 15-21 & 13   & 0.22   \\
13-23 & 18   & -0.41  &  & 15-21 & 18   & -0.22  \\
14-16 & 13   & 0.16   &  & 15-21 & 13   & 0.22   \\
14-16 & 15   & -0.16  &  & 15-21 & 18   & -0.22  \\
15-21 & 13   & 0.22   &  & 15-24 & 2    & 0.39   \\
15-21 & 18   & -0.22  &  & 15-24 & 3    & 0.18   \\
15-21 & 13   & 0.22   &  & 15-24 & 13   & 0.90    \\
15-21 & 18   & -0.22  &  & 15-24 & 18   & -1.47  \\
15-24 & 2    & 0.39   &  & 16-19 & 13   & 0.30    \\
15-24 & 3    & 0.18   &  & 16-19 & 16   & -0.30   \\ \cline{5-7}
15-24 & 13   & 0.90    &  &             &      &        \\
15-24 & 16   & -0.30   &  &             &      &        \\
15-24 & 21   & -1.17  &  &             &      &        \\
16-19 & 13   & 0.30    &  &             &      &        \\
16-19 & 16   & -0.30   &  &             &      &        \\ \cline{1-3}
\end{tabular}
\end{table}

\section{Conclusion}\label{sec:conclusion}
A novel methodology to solve an \gls{scopf} with short- and long-term post-contingency limits by using the \gls{imml} and Benders decomposition has been presented and applied to a case study. The case study suggests that the proposed methodology reduces the computational time significantly. 
The proposed methodology is specialized towards the \gls{scopf} framework, and is more complex to implement when compared to a standard direct solution method. This disadvantage needs to be balanced against the desired run time to solution. For system operation, reduced run time could mean closer to real-time analysis of critical contingencies if the system changes abruptly. And for planning, more scenarios can be explored within the same time-frame.

To further improve the run time of the proposed methodology, the iteration through contingencies could be executed in parallel. It should be noted that a parallel contingency overload search will often find more contingencies with overloads than if Benders' cuts are added in between, as some cuts mitigate overloads after other contingencies; solving again before adding cuts would mitigate this issue.

Further studies of the methodology's impact on system cost and reliability should be conducted. One useful addition to the \gls{scopf} is to include a risk index \cite{Wang2016} or chance constraints \cite{karangelos_iterative_2019}. A risk constraint can control the calculated risk in the system to balance the operational cost against risk averse operation. Also, extending the methodology to an AC \gls{scopf} would open the door to considering more types of operational limits.

\section{Appendix}\label{sec:Appendix}
The state of a power system can be analyzed using matrices. Here matrices are deduced for the linearized (DC) power flow equations. The first fundamental matrix is the connectivity matrix $\bm{\Phi}$ with dimensions ($B \times N$). This matrix shows how branches and nodes are connected. Specifically, it has a value of one for each node of a branch, while all other entries are zero.
\begin{equation}
    \Phi[l,n] = \ \left\{
    \begin{array}{@{}r@{\thinspace}l}
        1 ,&\quad\text{if } l_i = n \\
        -1 ,&\quad\text{if } l_j = n \\
        0,&\quad\text{otherwise}
    \end{array}
    \right. 
\end{equation}
where each branch $l$ has a node origin $l_i$ and end $l_j$, and $n$ is nodes in the set $\mathcal{N}$.

The second fundamental matrix is the diagonal susceptance matrix $\bm{\Psi}$ ($B \times B$), which contains the susceptance of all branches on its diagonal. All other values are zero.
\begin{gather}
    \Psi[k,k] = b_k \\
    \Psi[k,l] = 0 \quad\forall k \neq l
\end{gather}
where $b_k$ is the branch susceptance.
Combining these matrices results in the susceptance matrix $\bm{H}$.
\begin{equation}
    \bm{H} = \bm{\Phi}^\mathsf{T} \cdot \bm{\Psi} \cdot \bm{\Phi}
\end{equation}
where $\bm{\Phi}^\mathsf{T}$ is the transposed matrix of $\bm{\Phi}$.

The inverse susceptance matrix $\bm{X}$ is calculated from the susceptance matrix using Algorithm~\ref{alg:X}. $\bm{H}$ is singular when built from all system nodes, thus the reference node is eliminated before inverting $\bm{H}$. Further, $\bm{X}$ is used to solve the power flow and calculate the \gls{ptdf} matrix $\bm{\varphi}$.
\begin{algorithm}[tb]
    \DontPrintSemicolon
    \caption{Calculate the inverse susceptance matrix.}
    \label{alg:X}
    \KwIn{Susceptance matrix ($\bm{H}$) and reference node (ref)}
    \Begin{
    $\bm{H}[:,ref] \gets \bm{0}$\; 
    $\bm{H}[ref,:] \gets \bm{0}$\;
    $\bm{H}[ref,ref] \gets \bm{1}$\;
    $\bm{X} \gets inverse(\bm{H})$\;
    $\bm{X}[ref,ref] \gets 0$\;
    }
    \KwOut{Inverse susceptance matrix $\bm{X}$}
\end{algorithm}
\begin{gather}
    \bm{\theta} = \left[\bm{B}^*\right]^{-1} \cdot \bm{P} = \bm{X} \cdot \bm{P} \\ \label{eq:phi_norm} 
    \bm{\varphi} = \bm{\Psi} \cdot \bm{\Phi} \cdot \bm{X}
\end{gather}
where $\bm{\theta}$ is the node voltages angles, and $\bm{P}$ is the node injected active power. Branch power flow $\bm{F}$ can then be calculated two ways:
\begin{gather}\label{eq:flow_norm}
    \bm{F} = \bm{\Psi} \cdot \bm{\Phi} \cdot \bm{\theta} \\
    \bm{F} = \bm{\varphi} \cdot \bm{P}
\end{gather}

\bibliographystyle{IEEEtran}
\bibliography{references}
\end{document}